\documentclass{IOS-Book-Article}

\usepackage{mathptmx}
\usepackage{amsmath}
\usepackage{amssymb}
\usepackage{graphicx}
\usepackage{subfig}%
\usepackage{soul}\setuldepth{article}
\usepackage{comment}
\def\hb{\hbox to 11.5 cm{}}

\usepackage{amsthm}
\theoremstyle{plain}
\newtheorem{proposition}{Proposition}

\begin{document}

\pagestyle{headings}
\def\thepage{}
\begin{frontmatter}

\title{Tomography of nonlinear materials via the Monotonicity Principle}

\markboth{}{July 2023\hb}

\author[A]{\fnms{Vincenzo} \snm{Mottola}\orcid{0000-0002-8358-4544} %
\thanks{Corresponding Author: Vincenzo Mottola, vincenzo.mottola@unicas.it}},
\author[A]{\fnms{Antonio} \snm{Corbo Esposito}\orcid{0000-0003-4582-0128}}, %
\author[B]{\fnms{Gianpaolo} \snm{Piscitelli}\orcid{0000-0003-4591-5485}}
and
\author[A,C]{\fnms{Antonello} \snm{Tamburrino}\orcid{0000-0003-2462-6350}}%

\runningauthor{Antonello Tamburrino}
\address[A]{Dipartimento di Ingegneria Elettrica e dell'Informazione "M. Scarano", Università degli Studi di Cassino e del Lazio Meridionale, Italy}
\address[B]{Dipartimento di Scienze Economiche, Giuridiche, Informatiche e Motorie, Universit`a degli Studi di Napoli Parthenope, Italy}
\address[C]{Department of Electrical and Computer Engineering, Michigan State University, USA}

\begin{abstract}
In this paper we present a first non-iterative imaging method for nonlinear materials, based on Monotonicity Principle. Specifically, we deal with the inverse obstacle problem, where the aim is to retrieve a nonlinear anomaly embedded in linear known background.

The Monotonicity Principle (MP) is a general property for various class of PDEs, that has recently generalized to nonlinear elliptic PDEs. Basically, it states a monotone relation between the point-wise value of the unknown material property and the boundary measurements. It is at the foundation of a class of non-iterative imaging methods, characterized by a very low execution time that makes them ideal candidates for real-time applications.

In this work, we develop an inversion method that overcomes some of the peculiar difficulties in practical application of MP to imaging of nonlinear materials, preserving the feasibility for real-time applications. For the sake of clarity, we focus on a specific application, i.e. the Magnetostatic Permeability Tomography where the goal is retrieving the unknown (nonlinear) permeability by boundary measurements in DC operations. This choice is motivated by applications in the inspection of boxes and containers for security.

Reconstructions from simulated data prove the effectiveness of the presented method.
\end{abstract}

\end{frontmatter}
\markboth{July 2023\hb}{July 2023\hb}

\section{Introduction and statement of the problem}\label{sec:int}
Inverse problems and, in particular, electromagnetic tomography in the presence of nonlinear materials, is a quite new research theme~\cite{Lam20}. The general lack of contributions is ascribed to the complexity in dealing with such kind of problem. Indeed, retrieving the shape of an anomaly embedded in a known background material is a nonlinear and highly ill-posed problem, even in the presence of linear materials. When nonlinear material are considered, new challenges have to be properly addressed.

Here, we consider Magnetostatic Permeability Tomography (MPT), where the aim is to retrieve a nonlinear anomaly embedded in a linear known background material, throughout boundary measurements in (quasi)static operations. This technique has encountered great interest for the detection of metallic objects in boxes or containers, in security of ports and airports~\cite{Do18}. 

Let $\Omega$ be the region of tomographic inspection and $A\subset\Omega$ the region occupied by the nonlinear anomaly, the target magnetic permeability takes the form
\begin{displaymath}
    \mu_A(x,H)=\begin{cases}
        \mu_{NL}(H) & \text{in $A$} \\
        \mu_{BG} & \text{in $\Omega\setminus A$},
    \end{cases}
\end{displaymath}
where $\mu_{NL}$ is the nonlinear magnetic permeability characterizing the anomaly, while $\mu_{BG}$ is the linear and spatially uniform magnetic permeability related to the background.
If $\Omega$ is a simply connected domain and it is free from current densities, it is possible to introduce a magnetic scalar potential $\psi$, such that $\mathbf{H}=-\nabla\psi$. Then, the problem can be modelled as
\begin{equation}\label{eq:dirprob}
    \begin{cases}
        \nabla\cdot\left(\mu_A\left(x,|\nabla\psi_A(x)|\right)\nabla\psi_A(x)\right)=0 & \text{in $\Omega$} \\
        \psi_A(x)=f(x)                                                                 & \text{on $\partial\Omega$}.
    \end{cases}
\end{equation}
where $f\in X_{\diamond}=\{u\in H^{1/2}(\partial\Omega) | \int_{\partial\Omega}u(x)\,dx=0\}$ and $\mu_A$ has to meet some suitable, but very general, assumptions, in order to guarantee the existence and uniqueness of solutions. For a complete discussion we refer to~\cite{Co21, Co23}.

A key role in the framework of the inverse problem is played by the so called Dirichlet-to-Neumann (DtN) operator, mapping the imposed boundary potential to the measured normal component of the magnetic flux density on the boundary
\begin{displaymath}
    \Lambda_A:f\in X_{\diamond} \mapsto \mathbf{B}\cdot\hat{\mathbf{n}}|_{\partial\Omega}=-\mu_A\partial_n\psi_A|_{\partial\Omega}\in X_{\diamond}^{'},
\end{displaymath}
where $\hat{\mathbf{n}}$ is the outer unit normal defined on $\partial\Omega$ and $\partial_n$ the normal derivative in the outward direction. In this setting the excitation is given by a surface current density $\mathbf{J}_s$ on $\partial\Omega$. The magnetic imposed boundary potential is directly related to the surface current density, since it is possible to prove that $\mathbf{J}_s=-\hat{\mathbf{n}}\times\nabla_s \psi|_{\partial\Omega}$.

Recently, the Monotoncity Principle (MP) has been generalized to nonlinear materials, under very general assumptions~\cite{Co21, Co23}. MP states a monotone relation between the point-wise value of the unknown material property (the magnetic permeability for the case of interest) and the boundary measurements, expressed by a proper boundary operator. In the linear case, the MP is the foundation of a class of non-iterative imaging methods~\cite{Ta02} that exhibit excellent performances and feasibility in real-time applications, which is a feature as desirable as it is rare, in this field.

Let $T\subset\Omega$ an arbitrary region, termed as test anomaly, filled by the nonlinear magnetic permeability $\mu_{NL}$. Furthermore, let us assume that $\mu_{NL}(H)>\mu_{BG}(x)$, $\forall\,x\in\Omega$ and $\forall\,H>0$. From~\cite{Co21, Co23}, it is possible to state 
\begin{equation}\label{eq:mono1}
    T\subseteq A \Longrightarrow \overline{\Lambda}_T\leqslant\overline{\Lambda}_A,
\end{equation}
where $\overline{\Lambda}_T\leqslant\overline{\Lambda}_A$ is understood as $\langle \overline{\Lambda}_T(f),f\rangle\leq \langle \overline{\Lambda}_A(f),f\rangle \: \forall\,f\in X_{\diamond}$ and $\overline{\Lambda}_V$, with $V\subset\Omega$ is the average DtN operator, defined as
\begin{displaymath}
     \overline{\Lambda}_V:f\in X_{\diamond}\mapsto\int_0^1\Lambda_V(\alpha f)\,d\alpha\in X_{\diamond}^{'}.
\end{displaymath}
Proposition~\eqref{eq:mono1} is equivalent to 
\begin{equation}\label{eq:mono2}
    \overline{\Lambda}_T \not\leqslant \overline{\Lambda}_A \Longrightarrow T \not\subseteq A,
\end{equation}
where $\overline{\Lambda}_T \not\leqslant \overline{\Lambda}_A$ means that $\exists\,f\in X_{\diamond}\, : \,   \langle \overline{\Lambda}_T(f),f\rangle > \langle \overline{\Lambda}_A(f),f\rangle$.
Condition~\eqref{eq:mono2} is the foundation of the reconstruction method~\cite{Ta02}. In fact, it allows to infer if a test anomaly is included in $A$ or not, starting from boundary measurements. Hence, repeating the test in~\eqref{eq:mono2} for various test anomalies, it is possible to get an estimate $\tilde{A}$ of the anomaly $A$ as
\begin{equation}\label{eq:alg}
    \tilde{A}=\bigcup_{k}\{T_k \,|\, \overline{\Lambda}_{T_k}\leqslant\overline{\Lambda}_A\}.
\end{equation}

Despite the apparent simplicity of the algorithm in~\eqref{eq:alg}, there are several relevant challenges to be faced in the case of nonlinear materials. In the linear case, the average DtN operator is a linear one, and the monotonicity test in~\eqref{eq:mono2} can be easily performed by means of eigenvalues and eigenfunctions. In the nonlinear case, this is not the case, and there are no general mathematical tools to verify condition~\eqref{eq:mono2}. This represents one of the major barrier to the development of real-time imaging methods for nonlinear materials based on the MP.

In this paper, we propose initial numerical results for a new approach consisting in determining a set of \emph{optimal} boundary potential to verify~\eqref{eq:mono2}. These potentials are evaluated before the measurement process and by linear problems, guaranteeing the feasibility of the proposed method with real-time applications.

\section{Main idea and Proposed Method}
As mentioned in Section~\ref{sec:int}, the aim of this paper is to introduce an easy, systemic way to verify condition~\eqref{eq:mono2}. A possible solution is represented by performing the following nonlinear minimization problem
\begin{displaymath}
    \min_{f\in X_{\diamond}}\langle \overline{\Lambda}_A(f)-\overline{\Lambda}_T(f),f\rangle.
\end{displaymath}
where the goal is to determine if the minimum is negative or not. The minimization problem requires a huge number of iterations to converge, given the non-linearity and ill-posedness of the problem and, therefore, requires a huge number of measurements. For this reason suck kind of approach is certainly incompatible for real-time applications.

\begin{figure}[htp]
    \centering
    \subfloat[][]
    {\includegraphics[width=.35\textwidth]{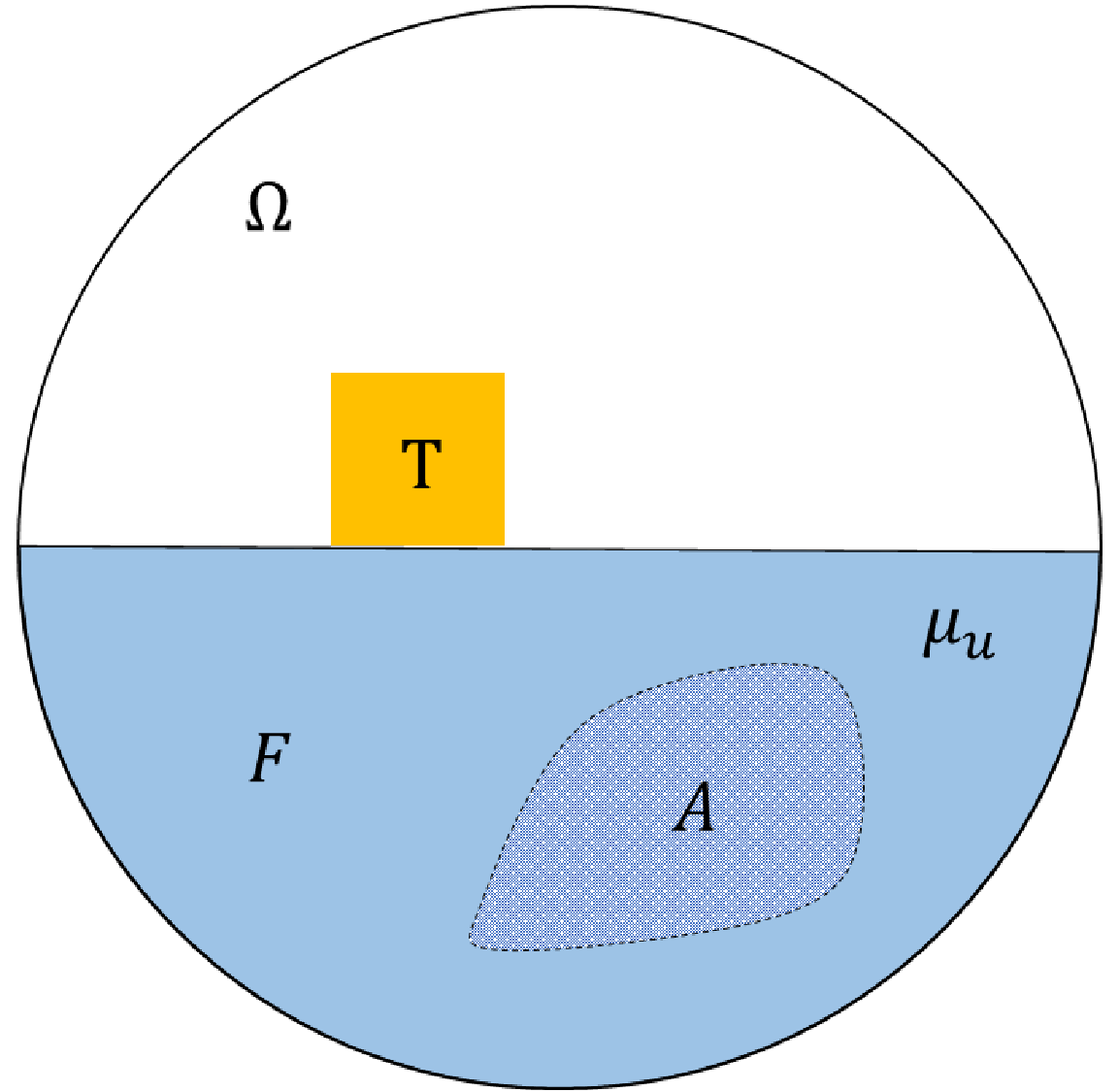}} \quad
    \subfloat[][]
    {\includegraphics[width=.45\textwidth]{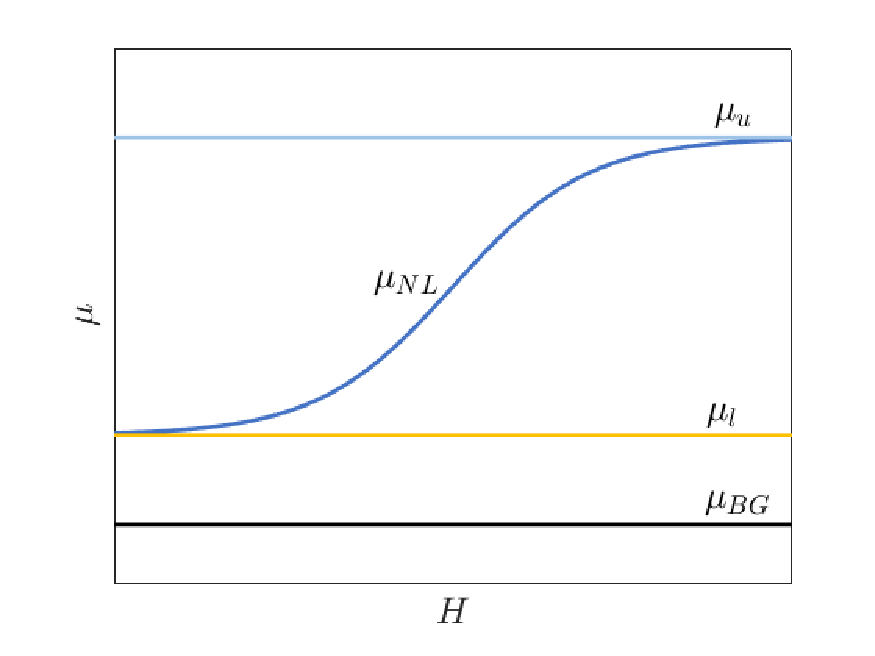}} \quad
    \caption{Left: Magnetic permeabilities involved. Right: Test anomaly $T$, Actual anomaly $A$ included in the region $F$.}
    \label{fig:perm}
\end{figure}

In this paper, we consider a different approach based on some physical observations. Let $\mu_u$ and $\mu_l$ two positive constants such that $\mu_l\leq\mu_{NL}(H)\leq\mu_u$, $\forall\,H>0$ (see Figure~\ref{fig:perm}). 

Given a certain test anomaly $T \not\subseteq A$, our aim is to find a proper boundary potential $f_0$ such that
\begin{displaymath}
    \langle \overline{\Lambda}_A(f_0)-\overline{\Lambda}_T(f_0),f_0\rangle<0\quad\text{if $T\not\subseteq A$}.
\end{displaymath}
Let $F\subset\Omega$ be an arbitrary region such that $A\subseteq F$, let $\mu_F^u$ be the magnetic permeability defined as
\begin{equation*}
\mu_F^u=\begin{cases}
    \mu^u & \text{in $F$}, \\
    \mu_{BG} & \text{in $\Omega\setminus F$}
\end{cases}    
\end{equation*}
and let  $\overline{\Lambda}_F^u$ be the related (linear) average DtN, let $\mu_T^l$ be the linear version of $\mu_T$ obtained when the test anomaly $T$ is filled by the linear material $\mu_l$ and let $\overline{\Lambda}_T^l$ be the related (linear) average DtN. It it is possible to prove the following result \cite{Mo23}.
\begin{proposition}
    Let $A\subset\Omega$ the anomalous region filled by the nonlinear magnetic permeability $\mu_{NL}$. Let $T\subset\Omega$ the region occupied by the test anomaly and let $F\subset\Omega$ such that $A\subseteq F$. If $\overline{\Lambda}_F^u-\overline{\Lambda}_T^l$ has negative eigenvalues, then any eigenfuctions $f_0$ related to negative eigenvalues ensures
    \begin{equation*}
        \langle \overline{\Lambda}_A(f_0)-\overline{\Lambda}_T(f_0),f_0\rangle<0.
    \end{equation*}
\end{proposition}

We are now in position to present our method for tomography of nonlinear materials:
\begin{enumerate}
    \item Given the geometry of $\Omega$, define a set of test anomalies $\{T_k\}_k$;
    \item For each test anomaly $T_k$, define a set of fictitious anomalies $\{F_{k,i}\}_i$ such that $\{F_{k,i}\}_i \cap T = \emptyset$; 
    \item for each test anomaly $T_k$ and for each fictitious anomaly $F_{k,i}$, compute and store the eigenfunction $f_{k,i}$ corresponding to the minimum (negative) eigenvalue of $\overline{\Lambda}_{F_{k,i}}^u-\overline{\Lambda}_{T,k}^l$;
    \item measure $\overline{\Lambda}_A$, for each $f_{k,i}$;
    \item the reconstruction of ${A}$ is given by $
        \tilde{A}=\bigcup \{T_k\,|\, \langle \overline{\Lambda}_A(f_{k,i})-\overline{\Lambda}_{T_k}(f_{k,i}),f_{k,i}\rangle\geq 0\: \forall\, i \}$.
\end{enumerate}

It is worth noting that only steps 4 and 5 need to be repeated at any new tomographic inspection. In other words, steps 1-3 can be carried out once for all, for a prescribed tomographic system.

\section{Numerical Examples}
In this section, we report some reconstructions obtained from numerically simulated data, corrupted by synthetic noise to avoid the inverse crime. Specifically, let $\overline{\Lambda}_{A,n}$ the noisy version of $\overline{\Lambda}_A$, we adopt the following noise model
\begin{displaymath}
    \langle\overline{\Lambda}_{A,n}(f),f\rangle=\langle\overline{\Lambda}_A(f),f\rangle+\eta\delta\xi,
\end{displaymath}
where $\eta$ is parameter controlling the noise level, $\xi \sim U(-1,1)$ is a uniformly distributed random variable, $\delta=\max_{f\in\mathbb{F}}\left\lvert \langle\overline{\Lambda}_A(f)-\overline{\Lambda}_0(f),f\rangle \right\rvert$ where $\mathbb{F}$ is the set of potential arising from step 3 of the method and $\overline{\Lambda}_0$ is the average DtN related to a configuration without defects (background only). In order to deal with noisy measurements we adopt the following reconstruction rule
    \begin{equation}\label{eqn:rn}
        \tilde{A}=\bigcup_k \{T_k\,|\, \langle \overline{\Lambda}_n(f_{k,i})-\overline{\Lambda}_{T_k}(f_{k,i}),f_{k,i}\rangle\geq -\eta\delta\: \forall\, i \}.
    \end{equation}
The~\eqref{eqn:rn} follows the approach proposed in~\cite{Ta02,Ha15,Ta16} for linear problems. It ensures, when $A$ can be described as union of same $T_k$, that $A\subseteq\tilde{A}$.

\begin{figure}[htp]
    \centering
    \subfloat[][Rectangle]
    {\includegraphics[width=.3\textwidth]{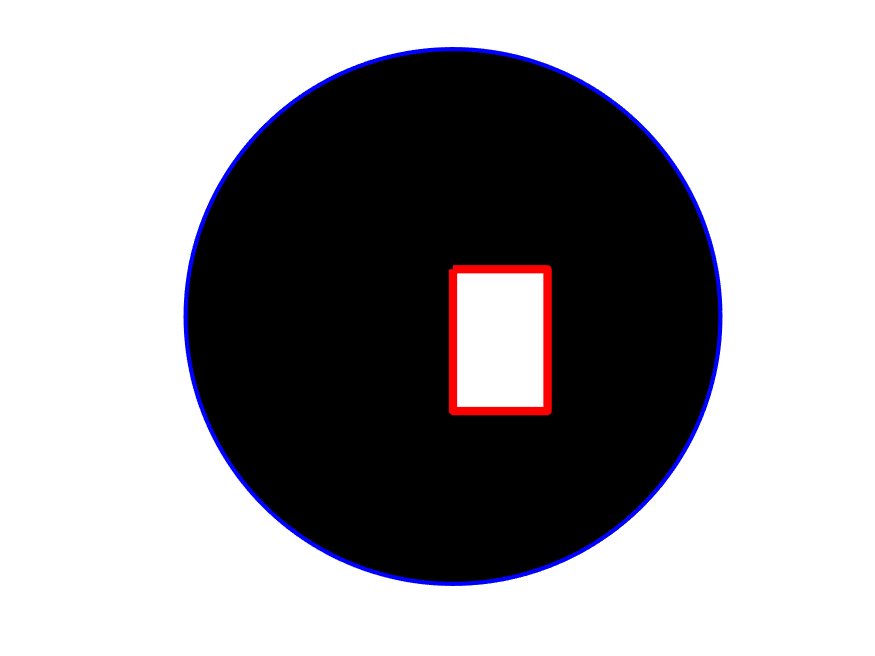}} \quad
    \subfloat[][L]
    {\includegraphics[width=.3\textwidth]{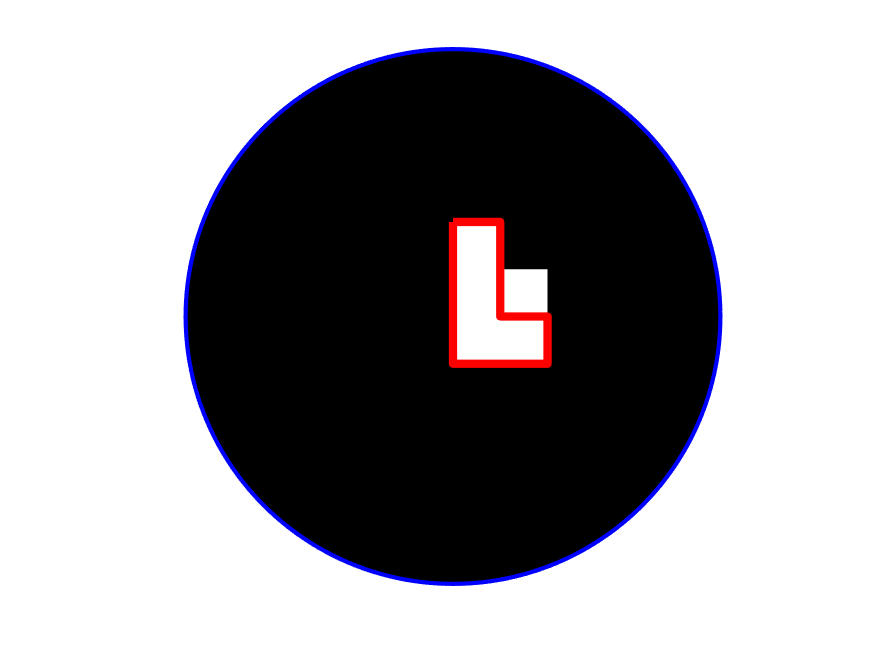}} \quad
    \subfloat[][C]
    {\includegraphics[width=.3\textwidth]{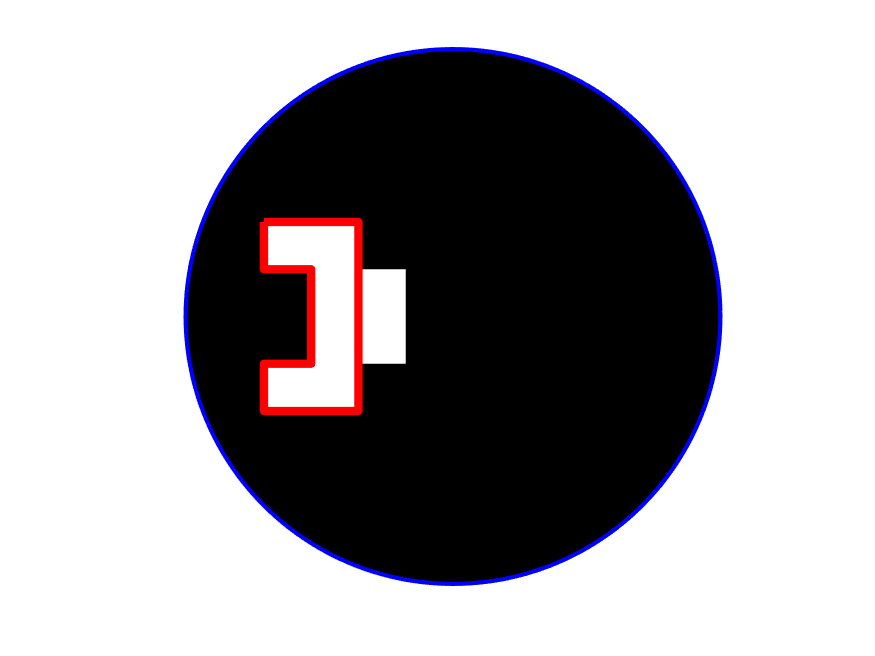}} \\
    \subfloat[][Circle]
    {\includegraphics[width=.3\textwidth]{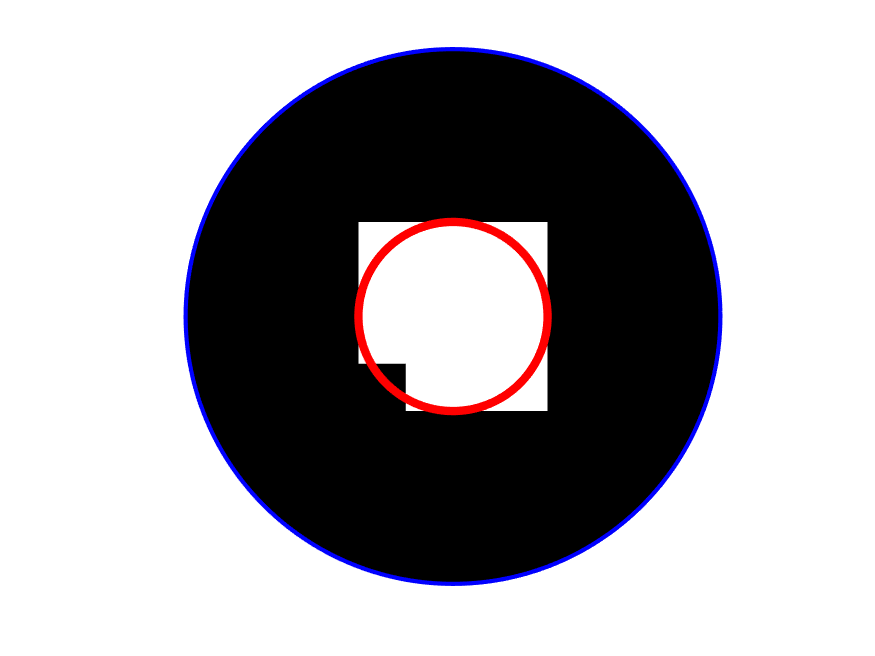}} \quad
    \subfloat[][Ellipse]
    {\includegraphics[width=.3\textwidth]{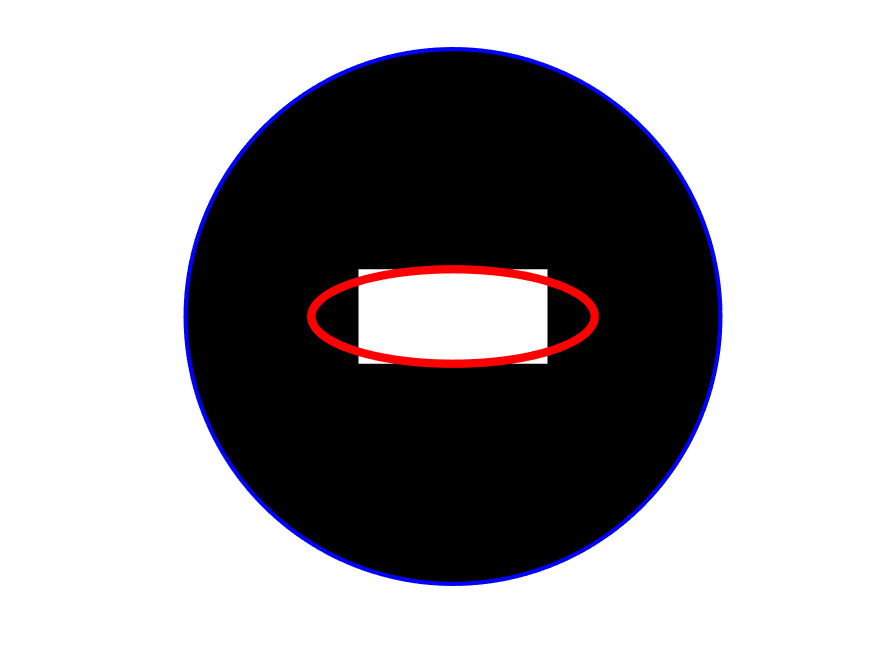}} \quad
    \subfloat[][Drop]
    {\includegraphics[width=.3\textwidth]{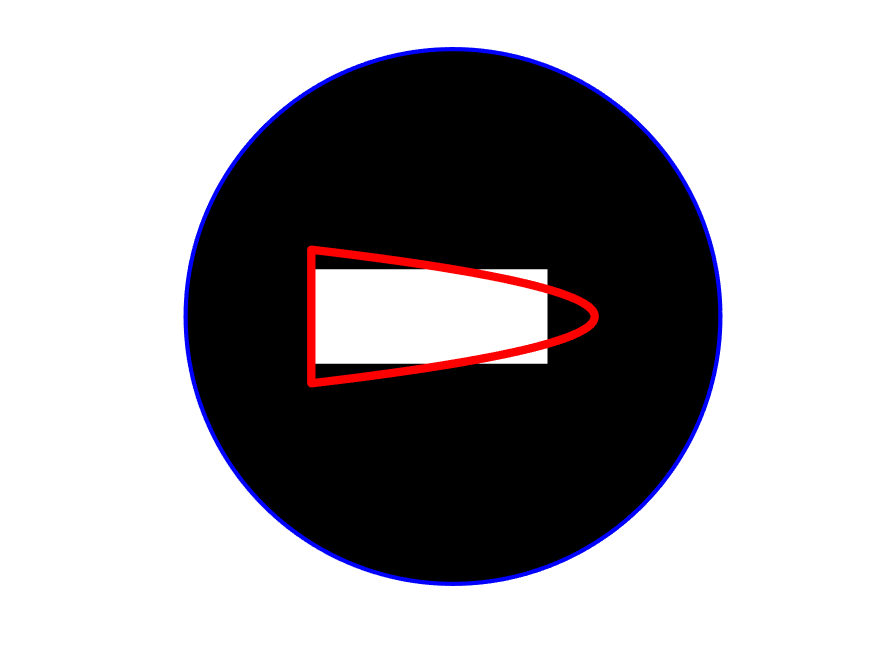}} \\
    \subfloat[][Bean]
    {\includegraphics[width=.3\textwidth]{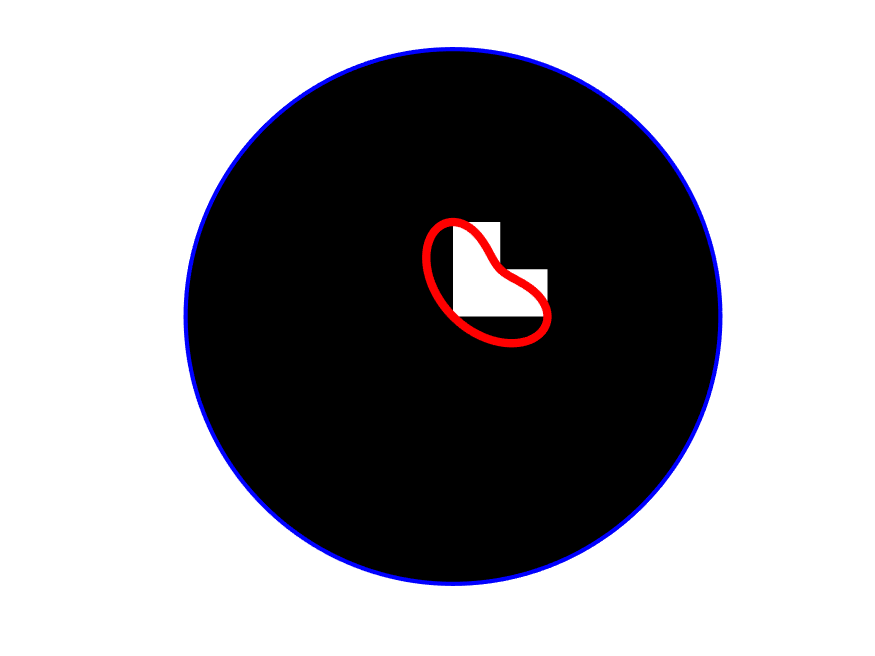}} \quad
    \subfloat[][Single test anomaly]
    {\includegraphics[width=.3\textwidth]{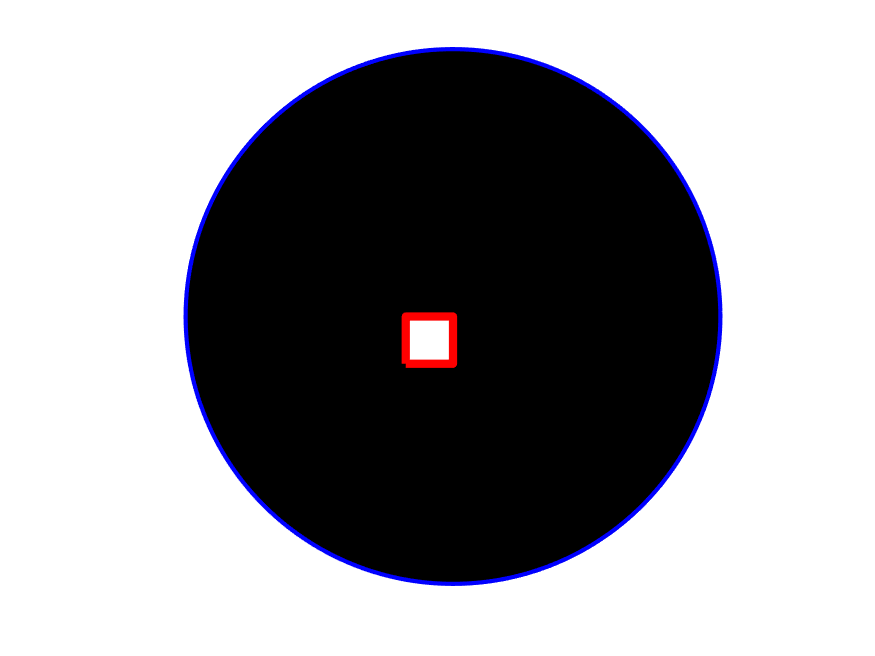}} \quad
    \subfloat[][Two connected components]
    {\includegraphics[width=.3\textwidth]{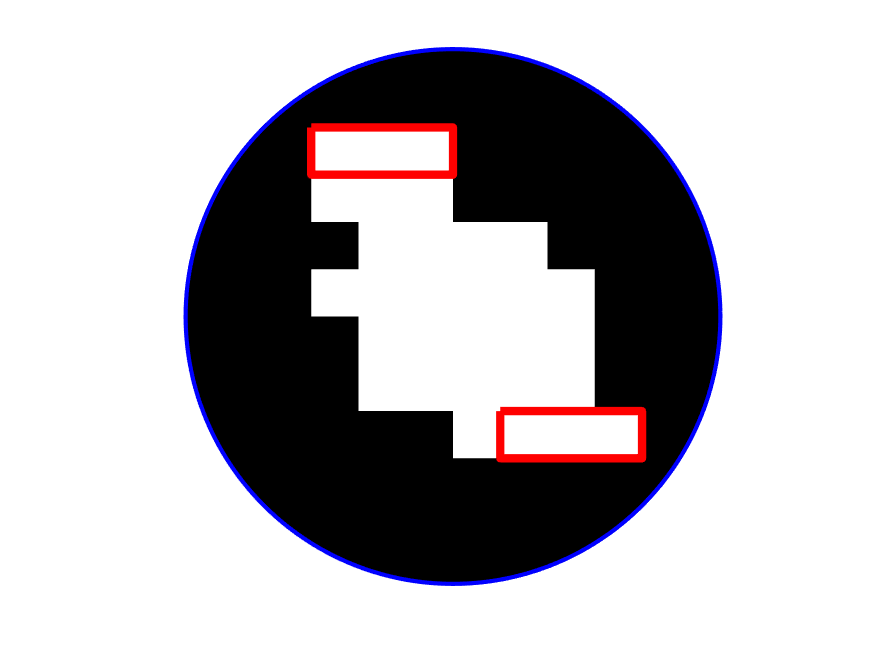}}
    \caption{Reconstructions carried out by the proposed method. In white the estimated anomaly $\tilde{A}$ while the boundary of the anomaly $A$ is in red. The noise level is equal to $\eta=10^{-3}$.}
    \label{fig:res}
\end{figure}

The domain $\Omega$ of tomographic inspection is represented by a circle of radius $30\,\text{cm}$. The aim is to retrieve the shape of ferromagnetic anomalies in air. In other words, $\mu_{BG}=\mu_0$, while the nonlinear anomaly $A$ is made by M330-50A Electrical Steel~\cite{Py13}. The problem is inspired by possible applications in the field of security in ports and airports, as reported in Section~\ref{sec:int} and reference therein. In Figure~\ref{fig:res}, reconstructions obtained for $\eta=10^{-4}$ are reported.

\section{Conclusions}
In this paper we present a non-iterative inversion method for tomography of nonlinear materials, based on Monotonicity Principle. We focus on Magnetostatic Permeability Tomography, where the aim is to retrieve the shape of a nonlinear anomaly embedded in a linear known background.
To the best of our knowledge, this is the first real-time inversion method specifically thought for the treatment of nonlinear materials.

The main challenge, in practical implementation of imaging methods for nonlinear materials based on Monotonicity Principle, is given by the nonlinearity of the average DtN, representing the boundary measurements. The nonlinearity causes a serious problem in finding the proper boundary potentials to carry put the Monotonicity Principle.

In this paper, we tested a new approach consisting in evaluating a set of proper boundary potential for the elementary monotonicity tests. These potentials allow to reveal if a test anomaly is not completely included in the actual anomaly. Moreover, they can be computed, before the measurements, by facing linear problems.

The strategy results in a method feasible for real-time applications. Numerical examples show the effectiveness of the method, even in presence of noise. The treatment of anomalies made by several connected components poses a relevant challenge and is currently under investigation.

\end{document}